\documentclass[smallextended]{svjour3} 
\usepackage{latexsym}
\usepackage{amsmath,amssymb,amsfonts}
\usepackage{tabu}
\usepackage{tikz}
\newif\ifstartedinmathmode
\newcommand\encircled[1]{%
  \relax\ifmmode\startedinmathmodetrue\else\startedinmathmodefalse\fi%
  \tikz[baseline,anchor=base]{%
  \node[draw,circle,outer sep=0pt,inner sep=.2ex]
    {\ifstartedinmathmode$#1$\else#1\fi};}%
}
\usepackage{fancybox}
\newcommand{\refe}[1]{(\ref{#1})}

\allowdisplaybreaks

\begin{document}

\title{Zeros of Jacobi and Ultraspherical  polynomials.\thanks{Dedicated to Richard (Dick) Askey, whose favorite polynomials were the Jacobi polynomials, in recognition of his seminal contributions to Special Functions and Orthogonal Polynomials}\thanks{The research of J. Arves\'u was funded by Agencia Estatal de Investigaci\'on of Spain, 
grant number PGC-2018-096504-B-C33. The research of K. Driver was funded by the National Research Foundation of South Africa, Grant Number 115332}
}


\titlerunning{Zeros of Jacobi polynomials}        

\author{J. Arves\'u   \and K. Driver \and  L.  Littlejohn}

\authorrunning{J. Arves\'u ,  K. Driver, L.  Littlejohn} 

\institute{J. Arves\'u \at
              Department of Mathematics, Universidad Carlos III de Madrid, 
Avda. de la Universidad, 30, 28911, Legan\'es, Spain
              \and
           K. Driver \at
              Department of Mathematics and Applied Mathematics, University of Cape Town, Cape Town 7708, South Africa
           \and 
           L.  Littlejohn \at  Department of Mathematics, Baylor University, One Bear Place 97328
Waco, TX 76798-7328, United States
}

\date{Received: date / Accepted: date}

\maketitle
\begin{abstract}Suppose $\{P_{n}^{(\alpha, \beta)}(x)\} _{n=0}^\infty $ is  a sequence of Jacobi polynomials with $ \alpha, \beta >-1.$  We discuss special cases of a question raised by Alan Sokal at OPSFA in 2019, namely, whether the zeros of  $ P_{n}^{(\alpha,\beta)}(x)$  and $ P_{n+k}^{(\alpha + t, \beta + s )}(x)$ are interlacing if  $s,t >0$ and $ k \in \mathbb{N}.$ We consider two cases of this question for Jacobi polynomials of consecutive degree and prove that the zeros of  $ P_{n}^{(\alpha,\beta)}(x)$  and $ P_{n+1}^{(\alpha, \beta + 1 )}(x),$   $ \alpha > -1, \beta > 0, $    $ n \in \mathbb{N},$  are partially, but in general not fully, interlacing depending on the values of $\alpha, \beta$ and $n.$    A similar result holds for the extent to which interlacing holds between the zeros of  $ P_{n}^{(\alpha,\beta)}(x)$  and $ P_{n+1}^{(\alpha + 1, \beta + 1 )}(x),$   $ \alpha >-1, \beta > -1.$  It is known that the zeros of the equal degree Jacobi polynomials $ P_{n}^{(\alpha,\beta)}(x)$  and $ P_{n}^{(\alpha - t, \beta + s )}(x)$ are interlacing for $ \alpha -t > -1, \beta > -1, $     $0 \leq t,s \leq 2.$ We prove that partial, but in general not full, interlacing of zeros holds between the zeros of  $ P_{n}^{(\alpha,\beta)}(x)$  and $ P_{n}^{(\alpha + 1, \beta + 1 )}(x),$  when  $ \alpha > -1, \beta > -1.$  We provide numerical examples that confirm that the results we prove cannot be strengthened in general. The  symmetric case $\alpha = \beta = \lambda -1/2$ of the Jacobi polynomials is also  considered. We prove that the zeros of the ultraspherical polynomials  $ C_{n}^{(\lambda)}(x)$  and $ C_{n + 1}^{(\lambda +1)}(x),$ $ \lambda >  -1/2$ are partially, but in general not fully, interlacing.  The interlacing of the zeros of the equal degree  ultraspherical polynomials  $ C_{n}^{(\lambda)}(x)$  and $ C_{n}^{(\lambda +3)}(x),$ $ \lambda >  -1/2,$ is also discussed.
\end{abstract}

\section{{Introduction}}  

\noindent If $\{p_{n}(x)\} _{n=0}^\infty $ is  a sequence of orthogonal polynomials with respect to the positive measure $d\mu(x)$ on a real interval $(a,b),$ it is well-known that the zeros of $p_n(x), n \geq 1,$ are real and simple and lie in $(a,b).$ Additionally, for $n \geq 1,$ the zeros of $p_n(x)$ and $p_{n+1}(x)$ interlace; that is to say between each pair of consecutive zeros of $p_{n+1}(x)$ there is one zero of $p_n(x).$ For example, the Jacobi polynomials $\{P_{n}^{(\alpha, \beta)}(x)\} _{n=0}^\infty,$ $ \alpha,\beta >-1,$  which are orthogonal on $(-1,1)$ with respect to the measure $d\mu_{\alpha,\beta}(x) = (1-x)^{\alpha} (1+x)^{\beta}dx,$  satisfy this  interlacing property. It is the subject of this paper to study new results on the zeros of Jacobi polynomials.

\medskip

\noindent The definition of interlacing of the zeros of two  polynomials of equal degree with real, simple zeros  is a natural extension of the definition of  interlacing of zeros of two polynomials of consecutive degree: If $x_{1} < x_{2} < \dots  < x_{n-1} < x_{n} $ and $y_{1} < y_{2} < \dots  < y_{n-1} < y_{n} $   denote the zeros of two polynomials of degree $n\geq 1,$ the zeros are said to interlace if either $x_{1} < y_{1} <  x_{2} < y_{2} < \dots  <  x_{n}< y_{n} $  or  $y_{1} < x_{1} <  y_{2} <  x_{2} < \dots  < y_{n} <  x_{n}.$ 

\medskip

\noindent In \cite{Ask},  Richard Askey proved that the zeros  of  $P_n^{(\alpha,\beta)}(x)$  and  $P_n^{(\alpha,\beta +1)}(x)$ are  interlacing for each $n \in \mathbb{N},$  $ \alpha > -1, \beta > -1,$   and  he  conjectured, based on graphical evidence, that  the zeros of  the Jacobi polynomials $P_n^{(\alpha,\beta)}(x)$ and  $P_n^{(\alpha,\beta +2)}(x)$   are  interlacing for each $n \in \mathbb{N}$, $ \alpha > -1, \beta > -1.$   Following Askey's result and his conjecture,  a more general result was proved in \cite{DrJoMb}, namely, that the zeros of  $P_n^{(\alpha,\beta)}(x)$  interlace with the zeros of  $P_n^{(\alpha -t,\beta +k)}(x)$ for each $n \in \mathbb{N}$ and each  $0 < t,k \leq 2$  provided  $ \alpha -t > -1, \beta > -1.$     In \cite{DIR}, Dimitrov, Ismail and Rafaeli consider interlacing properties of zeros of equal degree polynomials corresponding to different values of the parameter(s) within parameter-dependent orthogonal sequences by viewing parameter variations as perturbations of the weight function of orthogonality.

\medskip

\noindent A numerical example  given in \cite {DrJoMb} shows that when $\alpha = 1.26$ and $\beta = 1.85$,  the zeros of $P_4^{(\alpha,\beta)}(x)$  are $x_1 = -0.67979, \, x_2 =  -0.201233,\,  x_3 =  0.326414, \,  x_4 = 0.764756$ while the zeros of $P_4^{(\alpha +0.2,\beta +0.2)}(x)$ are $y_1 = -0.667543, \, y_2 = -0.197421 , \, y_3 = 0.317377, \,   y_4 =0.750436$ so full  interlacing between the two sets of zeros clearly breaks down.  However, in this example, $x_1 < y_1 < x_2 < y_2 $ while $ y_3 < x_3 < y_4 <x_4 $ which suggests that some, if not all,  of the zeros of $P_n^{(\alpha,\beta)}(x)$  may interlace with the zeros of  $P_n^{(\alpha + t,\beta +k)}(x)$  for every choice of  $n \in \mathbb{N},$  $ \alpha > -1, \beta > -1,$  and  $t = k = 1.$    

\medskip 

In Section \ref{section2}, we discuss the interlacing of the zeros of $P_n^{(\alpha,\beta)}(x)$  and  $P_{n+1}^{(\alpha,\beta +1)}(x),$  $ \alpha > -1, \beta > -1,$ and  the interlacing of the zeros of  $P_n^{(\alpha,\beta)}(x)$  and  $P_{n+1}^{(\alpha +1,\beta +1)}(x),$  $ \alpha > -1, \beta > -1.$  In Section \ref{section3}, we consider the interlacing of the zeros of  the two equal degree Jacobi polynomials  $P_n^{(\alpha,\beta)}(x)$  and $P_n^{(\alpha +1,\beta +1)}(x)$ where $ \alpha, \beta > -1.$  Section \ref{section4}  discusses the corresponding questions in the symmetric case $\alpha = \beta = \lambda -1/2$, namely, the interlacing of the zeros of ultraspherical  polynomials $C_n^{(\lambda)}(x)$  and $C_{n+1}^{(\lambda +1)}(x)$ as well as the interlacing of the zeros of  the equal degree ultraspherical polynomials $C_n^{(\lambda)}(x)$ and $C_n^{(\lambda +3)}(x),$ $\lambda > -1/2.$

\medskip

\section{ Zeros of  ${P_{n}^{(\alpha, \beta)}(x)}$  and ${P_{n+1}^{(\alpha +t, \beta +1)}(x)},$ $\alpha , \beta >-1,$ $t \in {\{0,1}\}.$ \label{section2}} 

Our first result considers the interlacing of  the zeros of  two Jacobi polynomials of consecutive degree $n$ and $n +1$ corresponding to the parameters $(\alpha,\beta)$  and $(\alpha, \beta +1)$ respectively.

\begin{theorem}\label{Theorem 2.1.} Suppose $\{P_{n}^{(\alpha, \beta)}(x)\} _{n=0}^\infty $ is  a sequence of Jacobi polynomials with $ \alpha > -1, \beta >0.$  Assume that $\alpha,\beta$ and $n$ are such that  $P_n^{(\alpha, \beta)}(x)$ and  $ P_{n+1}^{(\alpha, \beta +1)}(x)$ have no common zeros.  Let $ l_n = -1 + \dfrac {2(n+1)(\alpha +n +1)}{(\alpha + \beta +2n +2)(\alpha +\beta + 2n +3)}.$  If  $\{w_{i}\} _{i=1}^{n+1} $ are the zeros of  $ P_{n+1}^{(\alpha, \beta +1)}(x)$ in increasing order, $-1 < w_1 < w_2 < \dots < w_n < w_{n+1} < 1, $ each of  the $n+1$  open  intervals  $(-1, w_1),$ $(w_1, w_2),$ \dots, $(w_n, w_{n+1})$  contains either exactly  one simple zero of  $P_n^{(\alpha, \beta)}(x)$  or the point $ l_n $ but not both.  The zeros of $P_n^{(\alpha, \beta)}(x)$  and  $ P_{n+1}^{(\alpha, \beta +1)}(x),$ $\alpha > -1, \beta >0$   are  interlacing  if and only if the smallest (negative)  zero $w_1$  of  $ P_{n+1}^{(\alpha, \beta +1)}(x)$ satisfies  $w_1 >  l_n.$ 
\end{theorem}

\proof:  Eqn.(6) in \cite{DrJoMb} has a typographical error and should read:

\begin{multline} (\alpha + \beta + n +1)_2 (1+x) P_{n}^{(\alpha,\beta +2)}(x) = [ (\alpha + \beta + 2n +1)_2 (1+x)  - 2 n (\alpha +n)] P_{n}^{(\alpha,\beta)}(x) \\
 - (\alpha +n) [ (\alpha + 3\beta + 2n +4) + ( \alpha + \beta + 2n +2) x ]   P_{n-1}^{(\alpha,\beta +1)}(x), \label{(2.1)}
\end{multline}  
where $(r)_k = r(r+1)\dots (r+k-1)$, $(r)_{0}=1$, $k \in \mathbb{N},$ is the Pochhammer symbol. 

Replacing $n$ by $n+1$ and $\beta$ by $\beta -1$ in \refe{(2.1)} gives

\begin{multline} (\alpha + \beta + n +1)_2 (1+x) P_{n+1}^{(\alpha,\beta +1)}(x) = [ (\alpha + \beta + 2n +2)_2 (1+x)  - 2( n+1) (\alpha +n +1)] P_{n+1}^{(\alpha,\beta-1)}(x) \\
 - (\alpha +n +1) [ (\alpha + 3\beta + 2n +3) + ( \alpha + \beta + 2n +3) x ]   P_{n}^{(\alpha,\beta)}(x). \label{(2.2)}
\end{multline}

Evaluating \refe{(2.2)}  at successive zeros $w_k$ and $w_{k+1},$  $k \in {\{1,2,\dots,n}\},$  of  $ P_{n+1}^{(\alpha, \beta +1)}(x),$   $-1 < w_{1} <  w_{2} < \dots < w_{n +1} <1,$ we have 

\begin {equation}  a(w_k) a(w_{k+1})P_{n+1}^{(\alpha,\beta -1)}(w_k)P_{n+1}^{(\alpha,\beta -1)}(w_{k+1}) =  b(w_k) b(w_{k+1}) { k_n}^2   P_{n}^{(\alpha,\beta)}(w_k) P_{n}^{(\alpha,\beta)}(w_{k+1}), \label{(2.3)} \end{equation} 

where

\begin{align*} a(x) &= (\alpha + \beta + 2n +2) (\alpha + \beta + 2n +3)  (1+x)  - 2  (n+1)(\alpha +n+1),\\ 
b(x) &=  (\alpha + 3\beta + 2n +3) + ( \alpha + \beta + 2n +3) x ;\, \,  k_n = (\alpha +n +1).  \end{align*} 

We know from \cite[Th. 2.4]{DrJoMb} with $n$ replaced by $n+1$ and $\beta$ replaced by $\beta -1,$  that the zeros of ${P_{n+1}^{(\alpha, \beta +1)}(x)}$ and ${P_{n+1}^{(\alpha, \beta -1)}(x)}$ are interlacing provided $\alpha > -1$ and $\beta > 0.$
Therefore, $P_{n+1}^{(\alpha,\beta -1)}(w_k)P_{n+1}^{(\alpha,\beta -1)}(w_{k+1}) < 0$ for each $k \in {\{1,2,\dots ,n}\}.$  Also, $b(x)=0$ when $x = -1 -\frac{2\beta}{\alpha +\beta +2n +3}$ which is $ < -1$ for each $n \in \mathbb{N}, \alpha > -1, \beta >0.$ It follows that $b(x)$ has no zeros in $(-1,1)$ and $ b(w_k) b(w_{k+1})>0$ for each $k \in {\{1,2,\dots,n}\}.$ Clearly,  ${k_n}^2 >0.$   A calculation shows that $a(x) =0$ when  $x = l_n =  -1 + \dfrac {2(n+1)(\alpha +n +1)}{(\alpha + \beta +2n +2)(\alpha +\beta +2n +3)}$ which lies in $(-1,1).$  Therefore, $a(w_k) a(w_{k+1}) >0$  except for one value of  $k,$ $k \in {\{1,2,\dots ,n}\}.$  It follows that each interval $(w_k, w_{k+1}),$ $k \in {\{1,2,\dots,n}\},$ contains either $(i)$ exactly one zero of ${P_{n}^{(\alpha, \beta)}(x)}$ and not the point $l_n$   or $(ii)$ the point  $l_n$  and no zero of ${P_{n}^{(\alpha, \beta)}(x)}.$  Also,  a simple calculation confirms that   $l_n <0$ for each $n \in \mathbb{N}, \alpha >0, \beta >0.$ Therefore, each open interval  $(-1, w_1),(w_1, w_2), \dots, (w_n, w_{n+1})$  contains either exactly  $1$ simple zero of  $P_n^{(\alpha, \beta)}(x)$  or the point $l_n $ but not both.  This completes the proof.\qed 
 
\bigskip 

\begin{remark} 

Theorem \ref{Theorem 2.1.} can be equivalently reformulated as follows: Suppose  $\left\{P_{n}^{(\alpha, \beta)}(x)\right\}_{n=0}^\infty$ is  a sequence of Jacobi polynomials with  $\alpha > -1, \beta >0.$  Assume that $\alpha,\beta$ and $n$ are such that  $P_n^{(\alpha, \beta)}(x)$ and  $ P_{n+1}^{(\alpha, \beta +1)}(x)$ have no common zeros.  Let $l_n = -1 + \dfrac {2(n+1)(\alpha +n +1)}{(\alpha + \beta +2n +2)(\alpha +\beta + 2n +3)}.$  If  $\{w_{i}\}_{i=1}^{n+1}$ are the zeros of 
$P_{n+1}^{(\alpha, \beta +1)}(x)$ in increasing order, $-1 < w_1 < w_2 < \dots < w_n < w_{n+1} < 1,$  the zeros of  $(1+x)P_{n+1}^{(\alpha, \beta +1)}(x)$  and  $(x-l_n)P_n^{(\alpha, \beta)}(x)$  are interlacing. The zeros of $P_n^{(\alpha, \beta)}(x)$ and  $P_{n+1}^{(\alpha, \beta +1)}(x),$ $\alpha > -1, \beta >0$  are  interlacing  if and only if $l_n < w_1$, where $w_1$ is the smallest zero of  $P_{n+1}^{(\alpha, \beta +1)}(x).$  

\end{remark}

\bigskip

\begin{remark}

Our method of proof in Theorem 1 uses the interlacing property satisfied by the zeros of the Jacobi polynomials $ P_{n}^{(\alpha, \beta +1)}(x)$ and $P_{n}^{(\alpha, \beta -1)}(x)$ which is known to hold for each $n \in \mathbb{N}$,  $ \alpha > -1, \beta >0.$  However, the sequence $\{P_{n}^{(\alpha, \beta)}(x)\} _{n=0}^\infty $ is orthogonal for $ \alpha > -1, \beta > -1,$ so it is of interest to investigate whether interlacing holds between the zeros of $ P_{n}^{(\alpha, \beta +1)}(x)$ and $P_n^{(\alpha, \beta -1)}(x)$ when $ \alpha > -1$ and $ -1 <  \beta <0.$   From \cite[Th. 3 (ii)(b)] {BDR}, when $\alpha >-1, -1 <  \beta <0,$ the polynomial  $P_{n}^{(\alpha, \beta -1)}(x)$ is quasi-orthogonal of order $1$ and its $n$ zeros $\{y_{i}\} _{i=1}^{n} $ are real and distinct with 

\begin{equation}  y_1 < -1 < x_{1,n-1} < y_2 < x_{2,n-1} < y_3 < \dots < x_{n-1,n-1} < y_n <1. \label{(2.4)} \end{equation}  

where $\{x_{i, n-1}\} _{i=1}^{n-1} $ are the $n-1$ real, distinct zeros of $ P_{n-1}^{(\alpha, \beta)}(x)$ in $(-1,1).$

From \refe{(2.1)}  with $\beta$ replaced by $\beta -1,$

\begin{equation} h_n(1+x)  P_{n}^{(\alpha,\beta +1)}(x) = A(x) P_{n}^{(\alpha,\beta-1 )}(x)- B(x)(\alpha+n) P_{n -1}^{(\alpha,\beta)}(x), \label{(2.5)} \end{equation} 

where

\begin{equation*} h_n = (\alpha +\beta +n)_2,\end{equation*}

\begin{equation*} A(x) = (\alpha +\beta +2n)_2 (1+x) -2n(\alpha +n), \end{equation*}

\begin{equation*} B(x) = (\alpha +\beta +2n +1)x + (\alpha +3\beta +2n +1).\end{equation*}

Evaluating \refe{(2.5)} at successive zeros $y_k$ and $y_{k+1},$ $k \in {\{1,2,\dots,n-1}\}$  of  $P_{n}^{(\alpha,\beta-1 )}(x),$ we have

\begin{equation*} h_n^2 (1+y_k) (1+y_{k+1}) P_{n}^{(\alpha,\beta +1)}(y_k) P_{n}^{(\alpha,\beta +1)}(y_{k+1}) = B(y_k) B(y_{k+1}) (\alpha +n)^2 P_{n-1}^{(\alpha,\beta )}(y_k) P_{n -1}^{(\alpha,\beta)}(y_{k+1}). \end{equation*} 

From \refe{(2.4)}, $ (1+y_1) <0,$ $ (1+y_k) >0,$ for $k \in {\{2,\dots,n-1}\}.$  Also  $B(x) =0$ when $x = -1 -\dfrac{2\beta}{ (\alpha +\beta +2n +1) }$  which lies in the interval $(-1,1)$ when $-1<\beta <0.$ It follows that  partial but, in general, not full interlacing of zeros of  $ P_{n}^{(\alpha, \beta +1)}(x)$ and $P_n^{(\alpha, \beta -1)}(x)$ occurs when $ \alpha > -1$ and $ -1 <  \beta <0.$ 
\end{remark}

Table \ref{table-remark1-end} confirms that when $n=5$, $\alpha=10$, and $\beta=-0.1,$ the zeros of $P_n^{(\alpha, \beta +1)}(x)$ and  $P_n^{(\alpha, \beta -1)}(x)$ are interlacing.

\begin{table}[!htbp]
\caption{The zeros of $P_{n}^{(\alpha, \beta+1)}(x)$ and $P_{n}^{(\alpha, \beta-1)}(x)$ for $n=5, \alpha=10$, $\beta=-0.1.$ }
\begin{tabular}{|r|r|r|r|r|r|}
\hline
$P_{n}^{(\alpha,\beta+1)}(x)$ & $X_1$  & $X_2$  & $X_3$  & $X_4$  & $X_5$  \\
\hline
$P_{n}^{(\alpha,\beta-1)}(x)$ & $y_1$  & $y_2$  & $y_3$  & $y_4$  & $y_5$  \\
\hline\hline & & & & & \\ 
$P_{5}^{(10,0.9)}(x)$ & -0.9287 & -0.7588 & -0.5036 & -0.1807 & 0.1947 \\
$P_{5}^{(10,-1.1)}(x)$ & -1.0026 & -0.9112 & -0.6943 & -0.3704 & 0.0420 \\
\hline \end{tabular}
\label{table-remark1-end}
\end{table}

Table \ref{table-remark11-end} confirms that when $n=11$, $\alpha=1$, and $\beta=-1/2$  the zeros of $P_n^{(\alpha, \beta +1)}(x)$ and  $P_n^{(\alpha, \beta -1)}(x)$ are partially, but not fully, interlacing. The boxed points in Table \ref{table-remark11-end} highlight the interval with endpoints at a pair of consecutive zeros of $P_{n}^{(\alpha, \beta +1)}(x)$ that contains no zero of $P_{n}^{(\alpha, \beta-1)}(x).$

\begin{table}[!htbp]
\caption{The zeros of $P_{n}^{(\alpha, \beta+1)}(x)$ and $P_{n}^{(\alpha, \beta-1)}(x)$ for $n=11$ and $\alpha=1$, $\beta=-0.5$ .}
\begin{tabular}{|r|r|r|r|r|r|r|r|r|r|r|r|}
\hline
$P_{n}^{(\alpha,\beta+1)}(x)$ & $X_1$  & $X_2$  & $X_3$  & $X_4$  & $X_5$ & $X_6$ & $X_7$ & $X_8$ & $X_9$ 
& $X_{10}$ & $X_{11}$  \\
\hline
$P_{n}^{(\alpha,\beta-1)}(x)$ & $y_1$  & $y_2$  & $y_3$  & $y_4$  & $y_5$  & $y_6$ & $y_7$ & $y_8$ & $y_9$ 
& $y_{10}$ & $y_{11}$  \\
\hline\hline & & & & & & & & & & &\\ 
$P_{11}^{(1,0.5)}(x)$ & $\fbox{-0.967}$ & $\fbox{-0.871}$  & -0.718 & -0.518 & -0.284 & -0.031 & 0.224 & 0.464 & 0.674 & 0.840 & 0.951\\
$P_{11}^{(1,-1.5)}(x)$ & -1.005 & -0.969 & -0.855 & -0.675 & -0.443 & -0.177 & 0.103 & 0.374 & 0.617 & 0.811 & 0.942\\\hline \end{tabular}
\label{table-remark11-end}
\end{table}

\begin{corollary}\label{Corollary 2.1.}  Suppose $\{P_{n}^{(\alpha, \beta)}(x)\} _{n=0}^\infty $ is  a sequence of Jacobi polynomials with $ \alpha > 0, \beta >0$  and suppose that  $ P_{n+1}^{(\alpha +1, \beta)}(x)$ and $P_n^{(\alpha, \beta)}(x)$ have no common zeros.  If  $\{w_{i}\} _{i=1}^{n+1} $ are the zeros of  $ P_{n+1}^{(\alpha +1, \beta)}(x)$ in increasing order, $-1 < w_1 < w_2 < \dots < w_n < w_{n+1} < 1, $ each of  the $n+1$  open  intervals  $(w_1, w_2),$ \dots, $(w_n, w_{n+1}), (w_{n+1},1) $  contains either exactly  $1$ simple zero of  $P_n^{(\alpha, \beta)}(x)$  or the point $r_n = 1 - \dfrac {2(n+1)(\beta +n +1)}{(\alpha + \beta +2n +2)(\alpha +\beta + 2n +3)}$ but not both.  The zeros of $P_n^{(\alpha, \beta)}(x)$  and  $ P_{n+1}^{(\alpha +1, \beta)}(x),$ $\alpha > 0, \beta >0$   are  interlacing  if and only if the largest (positive)  zero $w_{n+1}$  of  $ P_{n+1}^{(\alpha, \beta +1)}(x)$ satisfies  $w_{n +1} <  r_n.$ 
\end{corollary}

\proof: The result follows from Theorem \ref{Theorem 2.1.} and the symmetry property  $P_n^{(\alpha, \beta)}(x) = (-1)^n   P_n^{(\beta, \alpha)}(-x),$  (see \cite[eqn.(4.1.3)] {Sze}).\qed

\begin{remark}\label{Remark 2.} 
\begin{itemize}

\item [(a)]  The point  $ l_n = -1 + \dfrac {2(n+1)(\alpha +n +1)}{(\alpha + \beta +2n +2)(\alpha +\beta + 2n +3)}$  depends on $n$, $\alpha$ and $\beta$ although the value of the parameter $\alpha$ is the same in the two polynomials  whose zeros we are comparing, namely, $ P_{n+1}^{(\alpha, \beta +1)}(x)$ and $ P_{n}^{(\alpha, \beta)}(x).$  
\item[(b)]  Each (fixed) value of the parameters $\alpha$ and $\beta$ with $\alpha > -1, \beta >0, $  generates an infinite sequence  $\{P_{n}^{(\alpha, \beta)}(x)\} _{n=0}^\infty $   of Jacobi polynomials that  is unique up to normalization. For $\alpha$ and $\beta$ fixed, $\alpha > -1, \beta >0, $  $l_n \rightarrow -1/2$ as $n \rightarrow \infty.$   It follows that breakdown of full interlacing of the zeros of $P_{n}^{(\alpha, \beta)}(x)$ and $ P_{n+1}^{(\alpha, \beta +1)}(x)$ is the most common, or the default, situation. For $n$ sufficiently large and any pair of fixed, finite values of $\alpha$ and $\beta,$ exactly $n-1$ of the $n$ zeros of $P_{n}^{(\alpha, \beta)}(x)$ interlace with the $n+1$ zeros of $ P_{n+1}^{(\alpha, \beta +1)}(x)$ and the point $l_n$ lies in the interval with endpoints at the pair of successive zeros of $ P_{n+1}^{(\alpha, \beta +1)}(x)$ that does not contain any zeros of $P_{n}^{(\alpha, \beta)}(x).$
\item[(c)] The result proved in Theorem \ref{Theorem 2.1.} is similar to the interlacing result for zeros of Laguerre polynomials $L_n^{(\alpha)}(x)$ and $L_{n+1}^{(\alpha + 1)}(x)$ proved in \cite[Th. 2.1]{ArDrLi}.
\end{itemize}
\end{remark}

Table \ref{Jacobi-1} and Table \ref{Jacobi-2} list the zeros of $P_{n}^{(\alpha, \beta)}(x)$ and $P_{n+1}^{(\alpha, \beta +1)}(x)$ for $n = 7$, $\alpha =6$, and $\beta \in {\{2,5,10,30,100}\}$. The numerical evidence confirms the results proved in Theorem \ref{Theorem 2.1.}. The boxed points in Table \ref{Jacobi-2} highlight the interval (if there is one) with endpoints at a pair of consecutive zeros of $P_{n +1}^{(\alpha, \beta +1)}(x)$ that contains the point  $l_n = -1 + \dfrac {2(n+1)(\alpha +n +1)}{(\alpha + \beta +2n +2)(\alpha +\beta + 2n +3)}$ and no zero of $P_{n}^{(\alpha, \beta)}(x).$ When $\beta$ is small compared with $n$, interlacing breaks down whereas, if $\beta >>n,$ interlacing holds.

\begin{table}[!htbp]
\caption{The zeros of $P_{n}^{(\alpha, \beta)}(x)$ for $n = 7$, $\alpha =6$, and $\beta \in {\{2,5,10, 30,100}\}$.}
\begin{tabular}{|r|r|r|r|r|r|r|}
\hline
 $x_1$  & $x_2$  & $x_3$  & $x_4$  & $x_5$  & $x_6$  & $x_7$ \\
\hline\hline & & & & & &\\
 -0.895  & -0.724  & -0.496  & -0.225  & 0.067  & 0.362  &  
   0.639  \\
 -0.765 & -0.553  & -0.310  & -0.048  & 0.218  & 0.473  &  0.706  \\
 -0.563 & -0.324  & -0.083  &  0.155 & 0.382  & 0.591  &  0.775  \\
 -0.036 &  0.180 & 0.363  & 0.523  & 0.663  & 0.783  &  0.883  \\
 0.537 & 0.652  & 0.740  &  0.811 & 0.870  & 0.918  & 0.956   \\
   \hline
\end{tabular}
\label{Jacobi-1}
\end{table}

\begin{table}[!htbp]
\caption{The zeros of $P_{n+1}^{(\alpha, \beta+1)}(x)$ for $n = 7$, $\alpha =6$, and $\beta \in {\{2,5,10,30,100}\}$.}
\begin{tabular}{|r|r|r|r|r|r|r|r|}
\hline
 $x_1$  & $x_2$  & $x_3$  & $x_4$  & $x_5$  & $x_6$  & $x_7$  & $x_8$ \\
\hline\hline & & & & & & &\\
-0.875  & $\fbox{-0.715}$  & $\fbox{-0.511}$  & -0.275  & -0.021  & 0.238  &
   0.486  & 0.713  \\
 $\fbox{-0.760}$ & $\fbox{-0.567}$  & -0.351  & -0.119  & 0.119  & 0.351  & 0.567 & 0.760   \\
 -0.580 & -0.363  & -0.145  & 0.072  & 0.282  & 0.479  & 0.657  & 0.812 \\
 -0.087 & 0.117  & 0.292  & 0.448  & 0.586  &  0.708 &  0.812 & 0.899 \\
 0.494 &  0.609 & 0.698  &  0.772 &  0.833 & 0.884  & 0.927  & 0.961 \\
 
   \hline
\end{tabular}
\label{Jacobi-2}
\end{table}

\begin{table}[!htbp]
\caption{The zeros  $x_1 < x_2 < \dots < x_{12}$  of  $P_{n}^{(\alpha, \beta)}(x)$ for $n = 12$, $\alpha =0.5$, and $\beta= 1$.}
\begin{tabular}{|r|r|r|r|r|r|r|r|r|r|r|r|}
\hline
 $x_1$  & $x_2$  & $x_3$  & $x_4$  & $x_5$  & $x_6$  & $x_7$  & $x_8$ & $x_9$ & $x_{10}$ & $x_{11}$ & $x_{12}$ 
\\
\hline\hline & & & & & & & & & & &\\
-0.958 &  -0.863 & -0.719 & -0.535 & -0.322 & -0.090 &
   0.147 &  0.375 & 0.583 & 0.757 & 0.890 & 0.972\\
     \hline
\end{tabular}
\label{Jacobi-3}
\end{table}

\begin{table}[!htbp]
\caption{The zeros $X_1 < X_2 < \dots < X_{13}$  of  $P_{n+1}^{(\alpha, \beta+1)}(x)$ for $n = 12$, $\alpha =0.5$, and $\beta= 1$.}
\begin{tabular}{|r|r|r|r|r|r|r|r|r|r|r|r|r|}
\hline
 $X_1$  & $X_2$  & $X_3$  & $X_4$  & $X_5$  & $X_6$  & $X_7$  & $X_8$ & $X_9$ & $X_{10}$ & $X_{11}$ & $X_{12}$ & $X_{13}$
\\
\hline\hline & & & & & & & & & & & &\\
 -0.940 &   -0.841 &  $\fbox{-0.705}$ & $\fbox{-0.537}$ &
  -0.345 & -0.138 & 
  0.076 & 0.286 &   0.482 & 0.657 & 0.802 & 0.910 & 0.977\\
  \hline
\end{tabular}
\label{Jacobi-4}
\end{table}

Table \ref{Jacobi-3} and Table \ref{Jacobi-4} list the zeros of  $P_{n}^{(\alpha, \beta)}(x)$ and $P_{n +1}^{(\alpha, \beta +1)}(x)$ for $n = 12$, 
$\alpha =0.5$, and $\beta= 1$. The interval with endpoints at successive zeros of $P_{n +1}^{(\alpha, \beta +1)}(x)$ where the interlacing of zeros breaks down is highlighted with boxes. Note that when $n=12,\alpha =0.5,\beta =1,$ we have $l_n= -0.552$ which lies in the boxed interval $(-0.705,-0.537)$. 

\medskip 

The next theorem considers the extent to which the zeros of two Jacobi polynomials are interlacing if the degree and each parameter value of the second polynomial is larger than the corresponding values of the degree and the parameters of the first polynomial. We introduce some notation in order to simplify the format of the equations used  to analyze the  interlacing of the zeros of  ${P_{n}^{(\alpha, \beta)}(x)}$  and ${P_{n+1}^{(\alpha+1, \beta +1)}(x)},$  $n \in \mathbb{N},$  $ \alpha > -1, \beta > -1.$   Let

\begin{align*} u(x) = u_{n,\alpha,\beta} (x) &: = \left(1-x^2\right) \frac{(\alpha +\beta +n +3)}{2}, \\ 
v = v_{n,\alpha,\beta} &:= \frac{2(\alpha +  n +2)(\beta +  n +2)}{(2n+ \alpha +\beta+4)}, \\
c(x) = c_{n,\alpha,\beta} (x) &:= x- \frac{(\alpha -\beta)}{(2n+ \alpha +\beta+4)},  \\
d = d_{n,\alpha,\beta} &:= 2(\alpha + \beta + n +2)(\alpha + \beta + 2n +2), \\
e = e_{n,\alpha,\beta} &:= (\alpha +\beta+ 2n +3), \\ 
f = f_{n,\alpha,\beta} &:=  (\alpha +\beta+ 2n +4)(\alpha +\beta+ 2n +2) = (e+1)(e-1) = e^2 -1, \\  
g = g_{n,\alpha,\beta} &:= 2(\alpha +  n +1)(\beta +  n +1)( e+1). \end{align*}

\begin{theorem}\label{Theorem 2.2.}  Suppose $\{P_{n}^{(\alpha, \beta)}(x)\} _{n=0}^\infty $ is  a sequence of Jacobi polynomials with $ \alpha > -1, \beta >-1.$   Assume that  $P_{n}^{(\alpha, \beta)}(x)$ and  $P_{n+1}^{(\alpha +1,\beta +1)}(x)$  have no common zeros. For $n\geq4,$ let $\{x_{k}\}_{k=1}^{n}$ denote the zeros of $P_{n}^{(\alpha,\beta)}(x)$ with $-1<x_{1}<x_{2}<\ldots<x_{n}<1.$ At least $n-3$ of the $n-1$ intervals with endpoints at successive zeros of $ P_{n}^{(\alpha, \beta)}(x)$  contain a simple zero of  $P_{n+1}^{(\alpha +1,\beta +1)}(x),$  $n \geq 4.$ The location of the $4$ remaining zeros of  $P_{n+1}^{(\alpha +1,\beta +1)}(x)$  is not possible to specify in general as different possibilities can occur depending on the values of $n, \alpha$, and $\beta.$ 
\end{theorem}

\proof: From  \cite [(3.15)]{DrJooJor} with $n$ replaced by $n+1$:

\begin{align} \left(1-x^2\right) \frac{(\alpha +\beta +n +3)} {2} P_{n+1}^{(\alpha +1,\beta +1)}(x) &=  \frac{2(\alpha +  n +2)(\beta +  n +2)}{(2n+ \alpha +\beta+4)} P_{n+1}^{(\alpha,\beta )}(x)  \notag\\
&- (n +2)  \left(x- \frac{(\alpha -\beta)}{ (2n+ \alpha +\beta+4) }\right)   P_{n +2}^{(\alpha,\beta)}(x).   
\label{(2.6)}\end{align} 

The three-term recurrence relation satisfied by Jacobi polynomials is

\begin{multline}  2(n +2)(n+ \alpha +\beta+2)( 2n+ \alpha +\beta+2)P_{n +2}^{(\alpha,\beta)}(x)   =  \\
 ( \alpha +\beta+2n +3) \left( \left(\alpha^2 -\beta^2\right) +x(\alpha +\beta +2n +4)(\alpha +\beta +2n +2)\right)  P_{n +1}^{(\alpha,\beta)}(x)  \\
-2 ( \alpha + n +1)(\beta +n +1)(\alpha +\beta +2n +4)  P_{n}^{(\alpha,\beta)}(x). \label{(2.7)}\end{multline}

Rewriting equations  \refe{(2.6)} and \refe{(2.7)} in terms of the notation introduced above gives

\begin{equation}   u(x)   P_{n +1}^{(\alpha +1,\beta +1)}(x) =  v   P_{n +1}^{(\alpha,\beta)}(x) - (n+2) c(x)  P_{n +2}^{(\alpha,\beta)}(x),  \label{(2.8)}\end{equation}  

and

\begin{equation}   (n+2) d   P_{n +2}^{(\alpha,\beta)}(x) =  e\left(\left(\alpha^2 -\beta^2\right) + (e^2-1)x\right)   P_{n +1}^{(\alpha,\beta)}(x) -  g  P_{n}^{(\alpha,\beta)}(x).  \label{(2.9)}\end{equation}

Multiplying \refe{(2.8)} by $d$ and \refe{(2.9)} by $-c(x)$ yields

\begin{equation}
d u(x)P_{n+1}^{(\alpha+1,\beta+1)}(x)- v d P_{n+1}^{(\alpha,\beta)}
(x)+dc(x)(n+2)P_{n+2}^{(\alpha,\beta)}(x)=0,\label{(2.10)}
\end{equation}
and
\begin{equation}
-dc(x)(n+2)P_{n+2}^{(\alpha,\beta)}(x)+ec(x)\left(  \alpha^{2}-\beta
^{2}+(e^{2}-1)x\right)  P_{n+1}^{(\alpha,\beta)}(x)-gc(x)P_{n}^{(\alpha
,\beta)}(x)=0.\label{(2.11)}
\end{equation}

Adding equations \refe{(2.10)} and \refe{(2.11)} gives
\begin{equation}
du(x)P_{n+1}^{(\alpha+1,\beta+1)}(x)+q(x)P_{n+1}^{(\alpha,\beta)}
(x)-gc(x)P_{n}^{(\alpha,\beta)}(x)=0,\label{(2.12)}
\end{equation}
where 
\begin{align}
q(x)  & :=ec(x)\left(  \alpha^{2}-\beta^{2}+(e^{2}-1)x\right)  -vd\nonumber\\
& =e(e^{2}-1)x^{2}+\left(  e(\alpha^{2}-\beta^{2})-e(\alpha-\beta
)(e-1)\right)  x-bd-\frac{e(\alpha-\beta)(\alpha^{2}-\beta^{2})}{e+1}
\label{(2.13)}\\
& =Ax^{2}+Bx+C.\nonumber
\end{align}

Evaluating \refe{(2.13)} at successive zeros $x_k$ and $x_{k+1},$  $k \in {\{1,2,\dots,n-1}\},$  of $ P_{n}^{(\alpha,\beta)}(x),$ we have
\begin{equation}
u(x_k) u(x_{k+1})  d^2   P_{n +1}^{(\alpha +1,\beta +1)}(x_k) P_{n +1}^{(\alpha +1,\beta +1)}(x_{k+1})=  q(x_k) q(x_{k+1}) P_{n +1}^{(\alpha,\beta)}(x_k)P_{n +1}^{(\alpha,\beta)}(x_{k+1}).\label{(2.14)}
\end{equation}

Since  $P_{n +1}^{(\alpha,\beta)}(x)$  and  $ P_{n}^{(\alpha, \beta)}(x)$ are polynomials of consecutive degree in an orthogonal sequence when $\alpha,\beta >-1$,  their zeros are simple and interlacing so that $P_{n+1}^{(\alpha
,\beta)}(x_{k})P_{n+1}^{(\alpha,\beta)}(x_{k+1})<0$ for $k \in {\{1,2,\dots,n}\}.$  Also, $u(x)>0$ for $x\in(-1,1)$ so $u(x_{k}) u(x_{k+1})>0$, while $q(x)$ is a quadratic function so that $q(x_{k})q(x_{k+1})\leq0$ for at most two values of $k\in\{1,2,\ldots,n-1\}.$ Hence,
from \refe{(2.14)}, $P_{n+1}^{(\alpha+1,\beta+1)}(x)$ has an odd number of zeros
in at least $n-3$ of the $n-1$ intervals $(x_{1},x_{2}),$ $(x_{2},x_{3}),\ldots,(x_{n-1},x_{n})$ when $n\geq4.$

We analyze the location of the zeros of $q(x).$ The discriminant $B^{2}-4AC$ of $q(x)$ is positive since
\begin{align}
B^{2}-4AC  & =e^{2}\left(  \alpha^{2}-\beta^{2}-(\alpha-\beta)(e-1)\right)
^{2}+4e(e^{2}-1)\left( vd+\frac{e(\alpha-\beta)\left(  \alpha^{2}-\beta
^{2}\right)  }{e+1}\right)  \nonumber\\
& =e^{2}(\alpha-\beta)^{2}\left(  \left(  \alpha+\beta-e+1\right)
^{2}+4(e-1)(\alpha+\beta)\right)  +4e(e^{2}-1)vd\nonumber\\
& =e^{2}(\alpha-\beta)^{2}(e-1+\alpha+\beta)^{2}+4e(e^{2}-1)vd>0,\label{(2.15)}
\end{align}
since $v,d>0$ and $e=\alpha+\beta+2n+3>2n+1.$ It follows that $q(x)$ has
two real roots $q_{\pm}:=q_{\pm}(n,\alpha,\beta).$ A tedious calculation shows that
\[
q(1)=4(n+1)(\beta+n+1)(\beta+n+2)>0,
\]
and
\[
q(-1)=4(n+1)(\alpha+n+1)(\alpha+n+2)>0.
\]
so the polynomial $q(x)$ is a U-shaped parabola with vertex $(x_{0},y_{0})$
where
\[
x_{0}=-\frac{B}{2A}=\frac{(\alpha-\beta)(n+1)}{(e^{2}-1)},
\]
and
\[
y_{0}=\frac{4AC-B^{2}}{4A}.
\]
From \refe{(2.15)}, we see that $y_{0}<0.$ If we can show that $x_{0}\in(-1,1),$
it will follow from the values of $q(x)$ at $x=\pm1$ and $y_{0}<0$ that the
real roots $q_{\pm}$ lie strictly between $-1$ and $1.$ We show that
$x_{0}<1;$ a similar proof shows that $x_{0}>-1.$ Since $\alpha,\beta>-1,$ we
see that
\begin{align*}
0  & <(\alpha+\beta)^{2}+4(n^{2}+n-1)\\
& <(\alpha+\beta)^{2}+5\alpha+7\beta+8+(3\alpha+5\beta+12)n+4n^{2},
\end{align*}
so that
\begin{align*}
e^2-1=(\alpha+\beta+2n+2)(\alpha+\beta+2n+4)  & =(\alpha+\beta)^{2}+6\alpha
+6\beta+8+(4\alpha+4\beta+12)n+4n^{2}\\
& >\alpha n+\alpha-\beta n-\beta=(\alpha-\beta)(n+1).
\end{align*}
It follows that $x_{0}<1.$ \qed

\begin{remark}\label{Theorem2_remark}
A natural question to ask is where the remaining four
zeros of $P_{n+1}^{(\alpha +1,\beta +1)}(x)$ lie. Since there may be `interference' from the two zeros of $q(x),$ it is difficult to exactly pinpoint where the other four zeros of $P_{n+1}^{(\alpha+1,\beta+1)}(x)$ lie
in $(-1,1)$ relative to the zeros of $P_{n}^{(\alpha,\beta)}(x)$. However, as illustrated by the examples
given in Table \ref{Jacobi-Theorem2}, there is an intimate connection between the location of these four zeros and the two roots of the quadratic polynomial $q(x)$. It appears that there are five possibilities for the location of these zeros:

\begin{itemize}

\item[(a)] each of the intervals $(x_{k},x_{k+1})$, $k\in \{1,2,\ldots
,n-1\},$ contains at least one root of $P_{n+1}^{(\alpha +1,\beta +1)}(x);$

\item[(b)] if there is a zero of $q(x)$ in $(-1,x_{1}),$ there is one zero
of $P_{n+1}^{(\alpha +1,\beta +1)}(x)$ in the interval $(-1,x_{1})$;

\item[(c)] if there is a zero of $q(x)$ in $(x_{n},1),$ there is one zero
of $P_{n+1}^{(\alpha +1,\beta +1)}(x)$ in the interval $(x_{n},1)$;

\item[(d)] if there is a zero of $q(x)$ in $(x_{k},x_{k+1})$ for some $k\in
\{1,2,\ldots ,n-1\},$ there are \underline{two} roots of $P_{n+1}^{(\alpha
+1,\beta +1)}(x)$ in $(x_{k},x_{k+1})$;

\item[(e)] it is possible that there is a root of $P_{n+1}^{(\alpha +1,\beta
+1)}(x)$ in $(-1,x_{1})$ but there is no root of $q(x)$ in this interval;
there is a similar situation for the interval $(x_{n},1)$ but we note that
it does not seem possible for both intervals to exclude roots of $q(x)$
but each contain a root of $P_{n+1}^{(\alpha +1,\beta +1)}(x).$ 
\end{itemize}

Regarding (c), if $q(x)$ has a zero in $(x_{k},x_{k+1}),$ we can show 
\begin{equation}
\mathrm{sgn}\left( P_{n+1}^{(\alpha +1,\beta +1)}(x_{k})\right) =\mathrm{sgn}%
\left( P_{n+1}^{(\alpha +1,\beta +1)}(x_{k+1})\right),   \label{(2.155)}
\end{equation}%
which implies that if $P_{n+1}^{(\alpha +1,\beta +1)}(x)$ has one root in $%
(x_{k},x_{k+1}),$ then it must have an even number of roots in this
interval. To show \refe{(2.155)}, we note that since%
\[
q(x_{k})q(x_{k+1})<0
\]%
and%
\[
P_{n+1}^{(\alpha ,\beta )}(x_{k})P_{n+1}^{(\alpha ,\beta )}(x_{k+1})<0,
\]%
we see from \refe{(2.13)} that%
\begin{align*}
\mathrm{sgn}\left( P_{n+1}^{(\alpha +1,\beta +1)}(x_{k})\right)  &=\mathrm{%
sgn}\left( du(x_{k})P_{n+1}^{(\alpha +1,\beta +1)}(x_{k})\right)  \\
&=\mathrm{sgn}\left( -q(x_{k})P_{n+1}^{(\alpha ,\beta )}(x_{k})\right)  \\
&=-\mathrm{sgn(}q(x_{k}))\mathrm{sgn}\left( P_{n+1}^{(\alpha ,\beta
)}(x_{k})\right)  \\
&=\mathrm{sgn}(q(x_{k+1}))\left( -\mathrm{sgn}\left( P_{n+1}^{(\alpha
,\beta )}(x_{k+1})\right) \right)  \\
&=\mathrm{sgn}\left( -q(x_{k+1})P_{n+1}^{(\alpha ,\beta )}(x_{k+1})\right) 
\\
&=\mathrm{sgn}\left( du(x_{k+1})P_{n+1}^{(\alpha +1,\beta
+1)}(x_{k+1})\right)  \\
&=\mathrm{sgn}\left( P_{n+1}^{(\alpha +1,\beta +1)}(x_{k+1})\right) .
\end{align*}
 Lastly, we note that,
for fixed $\alpha,\beta>-1,$ the zeros $q_{\pm}$ of $q(x)$ have the limiting
behavior
\[
\lim_{n\rightarrow\infty}q_{\pm}=\pm\frac{1}{\sqrt{2}}.
\]
\end{remark}

Table \ref{Jacobi-Theorem2} shows the zeros of the two polynomials $P_{n +1}^{(\alpha +1,\beta +1)}(x)$ and $ P_{n}^{(\alpha, \beta)}(x)$ for a selection of values $n$ and the parameters $(\alpha,\beta)$. Moreover, the two roots $q_{\pm}(n,\alpha,\beta)$ of the quadratic polynomial $q(x)$ are presented.

\medskip


\begin{table}[!htbp]
\caption{The zeros of $P_{n}^{(\alpha, \beta)}(x)$, $P_{n+1}^{(\alpha+1, \beta+1)}(x),$ and $q(x)$ for $n=8$, $\{(\alpha,\beta)\}=\{(-0.9,40), (-0.5,-0.5)\}$ and $n=6$, $(\alpha,\beta)=(27,29)$ are shown. }
\begin{tabu} to \textwidth { | X[c,m] | X[c,m] | X[c,m] | X[c,m] | X[c,m] | X[c,m] | X[c,m] | X[c,m] | X[c,m] | X[c,m] |}
\hline
$\hspace{-.1cm}P_{n+1}^{(\alpha+1, \beta+1)}$ & $X_1$  & $X_2$  & $X_3$ & $X_4$ & $X_5$ & $X_6$ & $X_7$  & $X_{8}$  & $X_{9}$   \\
$P_{n}^{(\alpha, \beta)}$ & $x_1$  & $x_2$  & $x_3$  & $x_4$ & $x_5$ & $x_6$ & $x_7$  & $x_8$  & $-$   \\
\hline
\end{tabu}
\begin{tabu} to \textwidth { | X[c,m] | X[c,m] | }
		\hline
		$q_{-}(n,\alpha,\beta)$ &
		$q_{+}(n,\alpha,\beta)$ \\
		\hline
	\end{tabu}
\begin{tabu} to \textwidth { | X[c,m] | X[c,m] | X[c,m] | X[c,m] | X[c,m] | X[c,m] | X[c,m] | X[c,m] | X[c,m] | X[c,m] |}\hline
$P_{9}^{(0.1, 41)}$  & 0.1776 & 0.3706 & 0.5263 & 0.6566 & 0.7650 & 0.8526 & 0.9198 & 0.9667 & 0.9932 \\
$P_{8}^{(-0.9, 40)}$ & 0.2810 & 0.4818 & 0.6384 & 0.7637 & 0.8615 & 0.9332 & 0.9792 & 0.9995 & $-$ \\
\hline
\end{tabu}
\begin{tabu} to \textwidth { | X[c,m] | X[c,m] | }
		\hline
		$q_{-}(8,-0.9, 40)=-0.9924$ &
		$q_{+}(8,-0.9, 40)=0.7742$ \\
		\hline
	\end{tabu}
\begin{tabu} to \textwidth { | X[c,m] | X[c,m] | X[c,m] | X[c,m] | X[c,m] | X[c,m] | X[c,m] | X[c,m] | X[c,m] | X[c,m] |}\hline
$P_{9}^{(\frac{1}{2}, \frac{1}{2})}$   & -0.9511 & -0.8090 & -0.5878 & -0.3090 & 0.0000 & 0.3090 & 0.5878 & 0.8090 & 0.9511 \\
$P_{8}^{(-\frac{1}{2}, -\frac{1}{2})}$ & -0.9808 & -0.8315 & -0.5556 & -0.1951 & 0.1951 & 0.5556 & 0.8315 & 0.9808 & $-$ \\
\hline
\end{tabu}
\begin{tabu} to \textwidth { | X[c,m] | X[c,m] | }
		\hline
		$q_{-}(8,-0.5, -0.5)=-0.7000$ &
		$q_{+}(8,-0.5, -0.5)=0.7500$ \\
		\hline
	\end{tabu}
\begin{tabu} to \textwidth { | X[c,m] | X[c,m] | X[c,m] | X[c,m] | X[c,m] | X[c,m] | X[c,m] | X[c,m] | X[c,m] | X[c,m] |}\hline	
$P_{7}^{(28, 30)}$  & -0.4142 & -0.2584 & -0.1128 & 0.0294 & 0.1709 & 0.3141 & 0.4654 & $-$ & $-$\\
$P_{6}^{(27, 29)}$  & -0.3751 & -0.2057 & -0.0470 & 0.1087 & 0.2654 & 0.4302 & $-$ & $-$ & $-$ \\	
\hline
\end{tabu}
\begin{tabu} to \textwidth { | X[c,m] | X[c,m] | }
		\hline
		$q_{-}(6,27, 29)=-0.9522$ &
		$q_{+}(6,27, 29)=0.9466$ \\
		\hline
	\end{tabu}	
	\label{Jacobi-Theorem2}
\end{table} 

\begin{remark}\label{Remark-2}  From Table 7, we see that the zeros $\{X_{i}\} _{i=1}^{9}$ of  $ P_{9}^{(0.1, 41)}(x),$  the zeros $\{x_{i}\} _{i=1}^{8}$ of  $ P_{8}^{(-0.9, 40)}(x)$ and $q_{+}=0.7742$ satisfy

\begin{equation*}X_1 <x_1 < X_2 <x_2 < X_3 <x_3 < X_4 <x_4 < X_5 <q_{+} < X_6<x_5 < X_7<x_6 < X_8<x_7 < X_9<x_8
\end{equation*}

In the symmetric case with $\alpha =\beta = -0.5$ and $n=8,$  the zeros $\{X_{i}\} _{i=1}^{9}$ of  $ P_{9}^{(0.5, 0.5)}(x),$  the zeros $\{x_{i}\} _{i=1}^{8}$ of  $ P_{8}^{(-0.5, -0.5)}(x)$ and $q_{+}=0.77$ and $q_{-}=-0.7$  satisfy

\begin{multline*}x_1 <X_1 < x_2 <X_2 <q_{-}< X_3 <x_3 < X_4 <x_4 < X_5(=0) <x_5 < X_6<x_6 \\
< X_7<q_{+} < X_8 < x_7 <X_9<x_8
\end{multline*}
Here, $q_{-}$ and $q_{+}$ play a role in the interlacing.

For $\alpha$ and $\beta$ large compared with $n,$ the zeros of ${P_{n}^{(\alpha, \beta)}(x)}$  and ${P_{n+1}^{(\alpha +1, \beta +1)}(x)},$ $\alpha , \beta >-1,$ are interlacing.

\end{remark}

\section{ Zeros of  ${P_{n}^{(\alpha, \beta)}(x)}$  and ${P_{n}^{(\alpha +1, \beta +1)}(x)},$ $\alpha , \beta >-1,$ $t \in {\{0,1}\}.$ \label{section3}} 

It was observed in \cite{DrJoMb} that numerical examples show that the zeros of ${P_{n}^{(\alpha, \beta)}(x)}$  and ${P_{n}^{(\alpha+1, \beta +1)}(x)}$ are not  interlacing for every $n \in \mathbb{N}$  and each choice of $ \alpha > -1, \beta > -1.$  However, a more precise statement can be made.

\medskip

\begin{theorem}\label{Theorem 3.1.} Suppose $\{P_{n}^{(\alpha, \beta)}(x)\} _{n=0}^\infty $ is  a sequence of Jacobi polynomials with $ \alpha > -1, \beta >-1.$  Assume that  $ P_{n}^{(\alpha, \beta)}(x)$ and  $P_{n}^{(\alpha +1,\beta +1)}(x)$  have no common zeros.  If $\{x_{i}\} _{i=1}^{n} $ are the zeros of  $ P_{n}^{(\alpha, \beta)}(x)$ in increasing order, $-1 < x_1 < x_2 < \dots < x_n < 1, $ at least $n-2$ zeros of ${P_{n}^{(\alpha+1, \beta +1)}(x)}$ lie strictly between a pair of consecutive zeros of  $ P_{n}^{(\alpha, \beta)}(x).$  The remaining two (simple) zeros of ${P_{n}^{(\alpha+1, \beta +1)}(x)}$ either both lie in one of the intervals with endpoints at a  pair of consecutive zeros of  $ P_{n}^{(\alpha, \beta)}(x)$ or one zero of  ${P_{n}^{(\alpha+1, \beta +1)}(x)}$  lies in the interval $(-1, x_1)$ and  one zero of  ${P_{n}^{(\alpha+1, \beta +1)}(x)}$  lies in the interval $( x_n, 1).$     
\end{theorem}

\proof: From  \cite [(3.15)] {DrJooJor}:

\begin{equation} 
J(x) P_{n}^{(\alpha +1,\beta +1)}(x) =  \frac{2(\alpha +  n +1)(\beta +  n +1)}{(\alpha +\beta+2n +2)} P_{n}^{(\alpha,\beta )}(x)  - H_n M(x)   P_{n +1}^{(\alpha,\beta)}(x),\label{(3.1)}
\end{equation} 

where

\begin{equation*} J(x) = \left(1-x^2 \right) \frac{\alpha + \beta + n +2 }{2};\quad H_n =  (n+1);\quad  M(x) = \left(x-\frac{\alpha -\beta}{ (\alpha +\beta +2 n+2)}\right).   \end{equation*} 

Evaluating \refe{(3.1)}  at successive zeros $ x_k$ and $ x_{k+1}$ of  $ P_{n}^{(\alpha, \beta)}(x)$  where $-1 < x_{1} <  x_{2} < \dots < x_{n} <1,$ we have, for $k \in \{1,2,\dots,n-1\},$

\begin{equation*}  J(x_k) J(x_{k+1})P_{n}^{(\alpha +1,\beta +1)}(x_k)P_{n}^{(\alpha +1,\beta +1)}(x_{k+1}) =  H_n^2 M(x_k) M(x_{k+1})  P_{n+1}^{(\alpha,\beta)}(x_k) P_{n +1}^{(\alpha,\beta)}(x_{k+1}).  \end{equation*}

Since  $P_{n +1}^{(\alpha,\beta)}(x)$  and  $ P_{n}^{(\alpha, \beta)}(x)$ are polynomials of consecutive degree in an orthogonal sequence, we know that the zeros of ${P_{n+1}^{(\alpha, \beta)}(x)}$ and ${P_{n}^{(\alpha, \beta)}(x)}$ are interlacing for $\alpha > -1$ and $\beta > -1$ so that $P_{n+1}^{(\alpha,\beta)}(w_k)P_{n+1}^{(\alpha,\beta)}(w_{k+1}) < 0$ for each $k \in \{1,2,\dots ,n-1\}.$  Also, ${H_n}^2 > 0$ for each $n \in \mathbb{N}$ while $M(x)$ is linear in $x$ so that $ M(x_k) M(x_{k+1})<0$ for at most one value of  $k \in {\{1,2,\dots,n-1}\}.$   Therefore,  $P_{n}^{(\alpha +1,\beta +1)}(x)$  has an odd number of  zeros, counting multiplicity,  in at least $n-2$ of the $n-1$ intervals with endpoints at successive zeros of $ P_{n}^{(\alpha, \beta)}(x).$   This accounts for at least $n-2$ zeros of  $P_{n}^{(\alpha +1,\beta +1)}(x)$  which has $n$ real simple zeros in the interval $(-1,1).$  The remaining (at most two) zeros of  $P_{n}^{(\alpha +1,\beta +1)}(x)$ must either both lie in the same interval with endpoints at successive zeros of  $ P_{n}^{(\alpha, \beta)}(x)$  or one zero of  ${P_{n}^{(\alpha+1, \beta +1)}(x)}$  lies in the interval $(-1, x_1)$ and  one zero of  ${P_{n}^{(\alpha+1, \beta +1)}(x)}$  lies in the interval $( x_n, 1).$   \qed

\medskip

\begin{remark}\label{Remark-3} It was proved in \cite[Th. 8]{DrJo1} that if $-2 < \alpha,\beta < -1,$ the Jacobi polynomial $P_{n}^{(\alpha,\beta)}(x)$ is quasi-orthogonal of order $2$ on $(-1,1)$ and the zeros of the (orthogonal) polynomial $P_{n}^{(\alpha +1,\beta +1)}(x)$ inerlace with the zeros of $ (x-\gamma) P_{n}^{(\alpha,\beta)}(x)$  provided $ P_{n}^{(\alpha,\beta)}(\gamma) \neq 0$ where $\gamma = \dfrac{\alpha - \beta}{\alpha + \beta +2n +2}.$ In the notation of Theorem \ref{Theorem 3.1.}, the zeros of $P_{n}^{(\alpha +1,\beta +1)}(x)$ interlace with the zeros of $ b(x) P_{n}^{(\alpha,\beta)}(x)$ provided $ P_{n}^{(\alpha,\beta)}(\gamma) \neq 0.$  However, the deduction made in \cite[Th.8]{DrJo1}, namely that the zeros of $P_{n}^{(\alpha +1,\beta +1)}(x)$ and  $ P_{n}^{(\alpha,\beta)}(x)$ are interlacing when $-2 < \alpha,\beta < -1$  cannot be made here because the reason interlacing holds for the "extra $2$" zeros is because  the  quasi-orthogonal order $2$ Jacobi polynomial $P_{n}^{(\alpha,\beta)}(x)$ has $2$ zeros lying outside the interval $(-1,1)$ whereas all the zeros of the (orthogonal) Jacobi polynomials $P_{n}^{(\alpha +1,\beta +1)}(x)$ and  $ P_{n}^{(\alpha,\beta)}(x)$ in Theorem \ref{Theorem 3.1.} lie in the interval $(-1,1).$ 
\end{remark}

\section{ Ultraspherical polynomials $\{C_{n}^{(\lambda)}(x)\} _{n=0}^\infty.$ \label{section4} } 

The one-parameter sequence of ultraspherical (also called Gegenbauer)  polynomials $\{C_{n}^{(\lambda)}(x)\} _{n=0}^{\infty} $ is the symmetric case of  the  sequence of Jacobi polynomials $\{P_{n}^{(\alpha, \beta)}(x)\} _{n=0}^\infty $ with $ \alpha= \beta = \lambda -1/2.$ For each $ \lambda > -1/2,$ the sequence  $\{C_{n}^{(\lambda)}(x)\} _{n=0}^\infty $  is  orthogonal on the interval $(-1,1)$  with respect to the (even) weight function $w(x) = {(1-x^2)}^{(\lambda -1/2)},$  the $n$ zeros of  $C_n^{(\lambda)}(x)$ are real, distinct, symmetric about the origin and lie in the open interval  $(-1,1)$ for each $n \in \mathbb{N}.$ The zeros of $C_{n-1}^{(\lambda)}(x)$ and $C_n^{(\lambda)}(x)$ are interlacing for each $n \in \mathbb{N}, n \geq 2.$

\medskip 

\noindent The  interlacing of zeros of two ultraspherical polynomials of equal or consecutive degree, $\lambda > -1/2,$ is  discussed in \cite {DrJo} and \cite {DIR}. Here, we consider interlacing of the zeros of  $C_{n}^{(\lambda)}(x)$ and  $C_{n+1}^{(\lambda +1)}(x).$ Note that an interlacing result can be deduced from Theorem 2 by putting $ \alpha= \beta = \lambda -1/2,$ but a stronger result can be proved directly. Note that when $n$ is odd, we consider interlacing of positive (respecively, negative) zeros of  two equal degree ultraspherical polynomials corresponding to different parameters since they have a common zero at the origin.  

\medskip

\begin{theorem}\label{Theorem 4.1.} Suppose  $\{C_{n}^{(\lambda)}(x)\} _{n=0}^\infty $  is a sequence of ultraspherical polynomials with $\lambda > -1/2.$  Assume that $ C_{n+1}^{(\lambda+1)}(x)$  and $ C_{n}^{(\lambda)}(x)$ have no common  zeros.  If  $\{w_{i}\} _{i=1}^{n+1} $ are the zeros of  $ C_{n+1}^{(\lambda+1)}(x)$ listed in increasing order, namely $-1 <  w_1 < w_2 < \dots < w_n < w_{n+1}  < 1,$   then each of  the $n +2$  open  intervals  $(-1, w_1), (w_1, w_2), \dots, (w_{n-1}, w_{n}), (w_{n},w_{n+1}), (w_{n+1},1)$ contains either exactly  one simple zero of  $C_{n}^{(\lambda)}(x)$  or the point $k_n$ or $-k_n$ where $ k_n = \sqrt{ \dfrac {(n+ 2\lambda +1)}{(2n + 2 \lambda + 2)}}$  but not both.  The zeros of $C_{n}^{(\lambda)}(x)$  and  $ C_{n+1}^{(\lambda +1)}(x),$ $\lambda > -1/2,$   are  interlacing  if and only if   $ -1 < -{k_n}< w_1$  or (equivalently)  $w_{n+1} < k_n < 1.$ 
\end{theorem}

\proof:  From \cite [(6)]{DrJo} with $\lambda$ replaced by $\lambda+1$:

\begin{equation}  (n +\lambda +1)  C_{n+1}^{(\lambda)}(x) =  \lambda\left( C_{n+1}^{(\lambda +1)}(x) -  C_{n-1}^{(\lambda +1)}(x)\right).  \label{(4.1)}\end{equation} 

From \cite [(7)]{DrJo}:

\begin{equation}  2\lambda (1-x^2)  C_{n-1}^{(\lambda +1)}(x) =  (2\lambda +n)x C_{n}^{(\lambda)}(x) - (n+1) C_{n+1}^{(\lambda)}(x).   \label{(4.2)}\end{equation} 

Multiplying \refe{(4.1)} by $2(1-x^2)$ and substituting from \refe{(4.2)}  gives  

\begin{equation} \left[ 2(1-x^2)(n +\lambda +1) - (n+1)\right]  C_{n+1}^{(\lambda)}(x) = 2\lambda (1-x^2) C_{n +1}^{(\lambda +1)}(x) - (2 \lambda +n) xC_{n}^{(\lambda)}(x).  \label{(4.3)} \end{equation} 

Assume that $n$ is even. Then $C_{n+1}^{(\lambda +1)}(x)$ and $C_{n+1}^{(\lambda)}(x)$  each have a simple zero at $x = 0$ and $n/2$ positive zeros. Evaluating \refe{(4.3)} at successive positive  zeros $z_k$ and $z_{k+1},$ 
$k \in \{ 1,2,\dots ,n/2\},$  of  $C_{n+1}^{(\lambda +1)}(x),$ we have 

\begin{multline*} \big[ 2(1-z_k^2)(n +\lambda +1) - (n+1)\big] [ 2(1-z_{k+1}^2)(n +\lambda +1) - (n+1)]  C_{n+1}^{(\lambda)}(z_k)  C_{n+1}^{(\lambda)}(z_{k+1})  \\
=  (2 \lambda +n)^2  z_k z_{k+1} C_{n}^{(\lambda)}(z_k) C_{n}^{(\lambda)}(z_{k+1}).   \end{multline*}

The positive (respectively, negative) zeros of  $ C_{n+1}^{(\lambda +1)}(x)$ and $ C_{n+1}^{(\lambda)}(x)$ are interlacing (\cite[Th.3.4]{DrJo}) so that $ C_{n+1}^{(\lambda)}(z_k)  C_{n+1}^{(\lambda)}(z_{k+1}) <0$ for each  
$k \in\{1,2,\dots ,n/2\},$  while $ (2 \lambda +n)^2 > 0. $  The quadratic equation $2(1-x^2)(n +\lambda +1) - (n+1) = 0$ has $2$ symmetric roots at $ x = k_n$ and $ x = -k_n$ where 
$k_n = \sqrt{ \dfrac{(n+ 2\lambda +1)}{(2n + 2 \lambda + 2)}}.$ The same argument as in Theorem \ref{Theorem 2.1.} completes the proof. \qed

\begin{remark}\label{Remark-Th4.2}  As $n \rightarrow \infty, k_n \rightarrow \frac{1}{\sqrt{2}}$ for $\lambda$ fixed and interlacing of the zeros of  $C_{n}^{(\lambda)}(x)$  and  $ C_{n+1}^{(\lambda +1)}(x),$ $\lambda > -1/2,$  does not hold for every $n \in \mathbb{N}.$  As $\lambda \rightarrow \infty,$ for $n$ fixed with $n << \lambda,$ $k_n \rightarrow 1,$ so we would expect that for $n$ small and $\lambda$ large, say $n=9$ and $\lambda = 4 000,$ interlacing will hold between all the zeros of  $C_{n}^{(\lambda)}(x)$  and  $ C_{n+1}^{(\lambda +1)}(x).$   \end{remark}

Tables \ref{Gegenbauer-Theorem4} and \ref{Gegenbauer-Theorem4-1} confirm the observations made. For $n$ small and $\lambda$  large, the zeros of $C_{n}^{(\lambda)}(x)$  and  $ C_{n+1}^{(\lambda +1)}(x)$ $\lambda > -1/2$   are  interlacing  but for $n$ (sufficiently) large and $\lambda$ fixed, full interlacing breaks down.

\begin{table}[!htbp]
\caption{The positive zeros of $C_{n}^{(\lambda)}(x)$ and $C_{n+1}^{(\lambda+1)}(x)$ for $n=9$ and $\lambda=4000$ are interlacing. }
\begin{tabular}{|r|r|r|r|r|r|}
\hline
$C_{n+1}^{(\lambda+1)}(x)$ & $X_1$  & $X_2$  & $X_3$  & $X_4$  & $X_5$   \\
\hline
$C_{n}^{(\lambda)}(x)$ & $x_1$  & $x_2$  & $x_3$  & $x_4$  & $-$   \\
\hline\hline & & & & &  \\
$C_{10}^{(4001)}(x)$   & 0.0054 & 0.0164 & 0.0278 & 0.0400 & 0.0543  \\   
$C_{9}^{(4000)}(x)$   & 0.0114 & 0.0232 & 0.0358 & 0.0504  &  $-$  \\
\hline \end{tabular}
\label{Gegenbauer-Theorem4}
\end{table}

\begin{table}[!htbp]
\caption{The positive zeros of $C_{n}^{(\lambda)}(x)$ and $C_{n+1}^{(\lambda+1)}(x)$ for $n=7$ and $\lambda=3.$  The   interval with endpoints at successive positive zeros of $C_{n+1}^{(\lambda+1)}(x)$ where the interlacing of zeros breaks down is highlighted with boxes.}
\begin{tabular}{|r|r|r|r|r|}
\hline
$C_{n+1}^{(\lambda+1)}(x)$ & $X_1$  & $X_2$ & $X_3$ & $X_4$\\
\hline
$C_{n}^{(\lambda)}(x)$ & $x_1$  & $x_2$ & $x_{3}$ & $-$\\
\hline\hline & & & & \\
$C^{(4)}_{8}(x)$  & 0.136   &   0.400  &  $\fbox{0.636}$  & $\fbox{0.830}$ \\   
$C^{(3)}_{7}(x)$   &  0.319   & 0.606 & 0.835 & $-$\\   
\hline \end{tabular}
\label{Gegenbauer-Theorem4-1}
\end{table}

\medskip

\noindent It is known (see \cite[Th.3.4]{DrJo}) that the positive (respectively, negative) zeros of the equal degree ultraspherical polynomials $C_{n}^{(\lambda)}(x)$  and  $ C_{n}^{(\lambda +2)}(x),$ $\lambda > -1/2$   are  interlacing for each $n \in \mathbb{N}, n \geq 2,$ and each $\lambda > -1/2.$ Numerical calculations confirm that the zeros of $C_{n}^{(\lambda)}(x)$  and  $ C_{n}^{(\lambda +3)}(x),$ $\lambda > -1/2,$  are not interlacing for every $n \in \mathbb{N}.$  Our next result proves that partial interlacing of zeros occurs.

\medskip

\begin{theorem}\label{Theorem 4.2.} Suppose  $\{C_{n}^{(\lambda)}(x)\} _{n=0}^\infty $  is a sequence of ultraspherical polynomials with $\lambda > -1/2,$ $\lambda \neq 0.$   Assume that  $C_{n}^{(\lambda +3)}(x)$  and $C_{n}^{(\lambda)}(x)$  have no common zeros other than the zero at the origin when $n$ is odd.  The zeros of  $C_{n}^{(\lambda +3)}(x)$  and $C_{n}^{(\lambda)}(x)$ satisfy a partial interlacing property.   
\end{theorem}

\proof:  From \cite[(9)]{DrJo} with $\lambda$ replaced by $\lambda+1$:

\begin{equation}  2(\lambda +2) (1-x^2) C_{n}^{(\lambda +3)}(x) =  A_{n,\lambda +1}(x) C_{n}^{(\lambda +1)}(x) + B_{n,\lambda +1}(x)   C_{n+1}^{(\lambda +1)}(x).  \label{(4.4)}\end{equation} 

where

\begin{multline}  A_{n,\lambda}(x) = \frac{( n + 2\lambda)}{ 2\lambda(1-x^2)}\left[(2\lambda +n +2) -(n+1)x^2\right],\quad   B_{n,\lambda}(x)= \frac{-( n + 1)x d_2^{(n,\lambda)}(x) }{ 2\lambda (1-x^2)}, \\
d_2^{(n,\lambda)}(x)= (2n+4\lambda +3) -2(n+\lambda +1)x^2. \label{(4.5)}\end{multline} 

The positive (respectively, negative) zeros of  $C_{n}^{(\lambda +1)}(x)$  and $C_{n}^{(\lambda +3)}(x) $ are interlacing  since their parameters differ by $2.$ (see \cite[Th.3.4]{DrJo}).  We look for an identity linking  $C_{n+1}^{(\lambda +1)}(x)$ with  $C_{n}^{(\lambda)}(x)$  and  $C_{n}^{(\lambda +1)}(x).$  From \refe{(4.5)}, we have

\begin{equation} B_{n,\lambda +1}(x)   C_{n+1}^{(\lambda +1)}(x) = \frac{-( n + 1)x}{ 2(\lambda+1) (1-x^2)} d_2^{(n,\lambda +1)}(x)  C_{n+1}^{(\lambda +1)}(x),  \label{(4.6)}\end{equation}

so it follows from \refe{(4.4)} and \refe{(4.6)} that

\begin{equation}  2(\lambda +2) (1-x^2) C_{n}^{(\lambda +3)}(x) =  A_{n,\lambda +1}(x)  C_{n}^{(\lambda +1)}(x)  - \frac{d_2^{(n,\lambda +1)}(x)}{ 2(\lambda+1) (1-x^2)}\left[(n+1) x  C_{n+1}^{(\lambda +1)}(x)\right].  \label{(4.7)}\end{equation} 

From \cite [(11)]{DrJo}

\begin{equation}  (n+1)x C_{n+1}^{(\lambda +1)}(x) = (2\lambda +n +2) C_{n}^{(\lambda+1)}(x) -2 (\lambda +1)(1-x^2) C_{n}^{(\lambda +2)}(x),   \label{(4.8)}\end{equation} 

and it follows from  \refe{(4.7)}, \refe{(4.8)}  that

\begin{multline*}  2(\lambda +2) (1-x^2) C_{n}^{(\lambda +3)}(x) =  A_{n,\lambda +1}(x)  C_{n}^{(\lambda +1)}(x) -\dfrac{d_2^{(n,\lambda +1)}(x)}{2(\lambda +1)(1-x^2)} \left((2\lambda +n+2)  C_{n}^{(\lambda +1)}(x)\right.\\
\left.  -2 (\lambda +1) (1-x^2) C_{n}^{(\lambda +2)}(x)\right).  \end{multline*} 

Therefore, 

\begin{equation*}  2(\lambda +2) (1-x^2) C_{n}^{(\lambda +3)}(x) = \left( A_{n,\lambda +1}(x) - \frac{(2\lambda +n +2) d_2^{(n,\lambda +1)}(x) }{2(\lambda +1)(1-x^2)}\right) C_{n}^{(\lambda +1)}(x)  +  d_2^{(n,\lambda +1)}(x) C_{n}^{(\lambda +2)}(x).  \end{equation*} 

A straightforward calculation using \refe{(4.5)} and \refe{(4.7)} shows that

\begin{equation*} A_{n,\lambda +1}(x) - \frac{(2\lambda +n +2) d_2^{(n,\lambda +1)}(x)}{2(\lambda +1)(1-x^2)} = \frac{- (n +2\lambda +2) (n +2\lambda +3)} { 2(\lambda +1)}, \end{equation*} 

and we obtain 

\begin{equation}  2(\lambda +2) (1-x^2) C_{n}^{(\lambda +3)}(x) =  d_2^{(n,\lambda +1)}(x)  C_{n}^{(\lambda +2)}(x)  -  h_n^{(\lambda +1)}  C_{n}^{(\lambda +1)}(x),  \label{(4.9)}\end{equation} 

where

\begin{equation} h_n^{(\lambda)} = \frac{ (n +2\lambda) (n +2\lambda +1)} { 2\lambda}. \end{equation}

Replacing $\lambda$ by $\lambda -1$ in \refe{(4.9)} gives

\begin{equation}  2(\lambda +1) (1-x^2) C_{n}^{(\lambda +2)}(x) =  d_2^{(n,\lambda)}(x)  C_{n}^{(\lambda +1)}(x)  -   h_n^{(\lambda)} C_{n}^{(\lambda)}(x).  \label{(4.11)}\end{equation}

Multiplying \refe{(4.11)} by $2(\lambda +1) (1-x^2):$

\begin{multline*} 4 (\lambda +1) (\lambda +2) {(1-x^2)}^2 C_{n}^{(\lambda +3)}(x) =  d_2^{(n,\lambda +1)}(x) 2(\lambda +1) (1-x^2) C_{n}^{(\lambda +2)}(x)  \\
-  2(\lambda +1) (1-x^2)  h_n^{(\lambda +1)}  C_{n}^{(\lambda +1)}(x),  \end{multline*}

or

\begin{multline}  4 (\lambda +1) (\lambda +2) {(1-x^2)}^2 C_{n}^{(\lambda +3)}(x) =  d_2^{(n,\lambda +1)}(x)\left[ d_2^{(n,\lambda)}(x) -  2(\lambda +1) (1-x^2)  h_n^{(\lambda +1)}\right] C_{n}^{(\lambda +1)}(x) \notag\\
- h_n^{(\lambda)} d_2^{(n,\lambda +1)}(x) C_{n}^{(\lambda)}(x).   \end{multline}  

Therefore  

\begin{equation*}  4 (\lambda +1) (\lambda +2) {(1-x^2)}^2 C_{n}^{(\lambda +3)}(x) =  H_4^{(n,\lambda)}(x) C_{n}^{(\lambda +1)}(x) - h_n^{(\lambda)} d_2^{(n,\lambda +1)}(x) C_{n}^{(\lambda)}(x),  
\end{equation*}  

where

\begin{equation} H_4^{(n,\lambda)}(x) = \left[ d_2^{(n,\lambda +1)}(x) d_2^{(n,\lambda)}(x) 
-  2(\lambda +1) (1-x^2)  h_n^{(\lambda +1)}\right].\label{(4.12)}\end{equation}

Evaluating \refe{(4.12)} at successive positive (respectively, negative) zeros $x_k$ and $x_{k+1}$ of $C_{n}^{(\lambda)}(x),$ we know that $C_{n}^{(\lambda +1)}(x_k)C_{n}^{(\lambda +1)}(x_{k+1}) <0$ for each $k \in {\{1,2,\dots n-1}\}$. Since $H_4^{(n,\lambda)}(x)$ is an even polynomial of degree $4$ in $x,$  it follows that interlacing of zeros may break down between the zeros of  $C_{n}^{(\lambda +3)}(x)$  and $C_{n}^{(\lambda)}(x)$  in at most two intervals with endpoints at successive positive (respectively, negative) zeros of $C_{n}^{(\lambda)}(x).$ \qed

Numerical evidence shows that the zeros of $C_{n}^{(\lambda +3)}(x)$  and $C_{n}^{(\lambda)}(x)$ interlace when $n$ is small and $\lambda$ is large (see Table \ref{Gegenbauer-final}). However, the interlacing of zeros
breaks down when $n$ is large relative to $\lambda$ (see Table \ref{Gegenbauer-final2}).

\begin{table}[!htbp]
\caption{The positive zeros of $C_{n}^{(\lambda+3)}(x)$ and $C_{n}^{(\lambda)}(x)$ for $n=11$ and $\lambda=100.$ }
\begin{tabular}{|r|r|r|r|r|r|}
\hline
$C_{n}^{(\lambda+3)}(x)$ & $X_1$  & $X_2$ & $X_3$ & $X_4$ & $X_{5}$\\
\hline
$C_{n}^{(\lambda)}(x)$ & $x_1$  & $x_2$   & $x_3$ & $x_4$ & $x_5$\\
\hline\hline & &  & &  & \\
$C_{11}^{(103)}(x)$     & 0.0631 & 0.1270 & 0.1929 & 0.2628 & 0.3420  \\
$C_{11}^{(100)}(x)$     & 0.0640 & 0.1288 & 0.1956 & 0.2664 & 0.3465 \\   
\hline \end{tabular}
\label{Gegenbauer-final}
\end{table}

\begin{table}[!htbp]
\caption{The positive zeros of $C_{n}^{(\lambda+3)}(x)$ and $C_{n}^{(\lambda)}(x)$ for $n=15$ and $\lambda=-1/4.$ The interval with endpoints at successive zeros of $C_{n}^{(\lambda+3)}(x)$ where the interlacing of zeros
breaks down is highlighted with boxes.}
\begin{tabular}{|r|r|r|r|r|r|r|r|}
\hline
$C_{n}^{(\lambda+3)}(x)$ & $X_1$  & $X_2$  & $X_3$  & $X_4$  & $X_5$ & $X_6$ & $X_7$ \\
\hline
$C_{n}^{(\lambda)}(x)$   & $x_1$  & $x_2$  & $x_3$  & $x_4$  & $x_5$ & $x_6$ & $x_7$ \\
\hline\hline & & & & & & & \\
$C_{15}^{(11/4)}(x)$   & 0.177   &  0.349   &  0.510   & 0.655    & 0.780   & $\fbox{0.880}$   & $\fbox{0.953}$ \\   
$C_{15}^{(-1/4)}(x)$   & 0.212   &  0.414   &  0.597   & 0.753    & 0.875   & 0.958   & 0.997\\
\hline \end{tabular}
\label{Gegenbauer-final2}
\end{table}

\begin{remark}   If we consider $n \in \mathbb{N}$ as the (fixed) parameter, equation \refe{(4.11)} is a three-term recurrence relation in  $\lambda.$  It is interesting to observe that the so-called Favard coefficient is equal to $h_n^{(\lambda)}$ which is positive for $\lambda > -1/2,$  $n \geq 2.$  
\end{remark}

\begin{acknowledgements}
Jorge Arves\'u and Kathy Driver wish to thank the Mathematics Department at Baylor University for hosting their visits in Fall 2019 which stimulated this research.
\end{acknowledgements}

\end{document}